\documentclass[11pt, leqno]{amsart}

\usepackage{amssymb}
\usepackage{amsmath}
\usepackage{amsthm}

\newtheorem{thm}{Theorem}[section]
\newtheorem{lemma}[thm]{Lemma}
\newtheorem{prop}[thm]{Proposition}
\newtheorem{cor}[thm]{Corollary}

\newcommand{\bC}{\mathbb{C}}
\newcommand{\bR}{\mathbb{R}}
\newcommand{\bZ}{\mathbb{Z}}
\newcommand{\bP}{\mathbb{P}}
\newcommand{\bN}{\mathbb{N}}
\newcommand{\bQ}{\mathbb{Q}}

\newcommand{\cF}{\mathcal{F}}
\newcommand{\cL}{\mathcal{P}}
\newcommand{\co}{\mathcal{O}}
\newcommand{\cH}{\mathcal{H}}

\newcommand{\cZ}{\mathcal{Z}}

\newcommand{\cP}{\Lambda}
\newcommand{\lb}{\mathcal{L}}
\newcommand{\cN}{N}
\newcommand{\fH}{\mathfrak{H}}

\newcommand{\fU}{\mathfrak{U}}

\newcommand{\Hom}{\mathrm{Hom}}
\newcommand{\Pic}{\mathrm{Pic}}

\newcommand{\ord}{\mathrm{ord}}
\newcommand{\Emb}{\mathrm{Emb}}
\newcommand{\R}{\mathrm{R}(\mathfrak{H}_{\varphi})}

\newcommand{\tp}{\tilde{\varphi}}

\newcommand{\ba}{\bar\partial}

\newcommand{\pgl}{PGL(2,\mathbb{C})}
\newcommand{\lp}{\Lambda_{\varphi}}
\newcommand{\op}{\mathrm{ord}(\varphi)}

\newcommand{\h}{\Hom(L\bC^{\ast},\mathbb{C}^{\ast})}
\newcommand{\hi}{\mathrm{Pic}(L\bP_1)^{LPGL(2,\bC)}}
\newcommand{\hp}{\mathfrak{H}_{\varphi}}
\newcommand{\hs}{\mathfrak{H}_{\varphi}(\sigma)}
\newcommand{\s}{H^0(L\mathbb{P}_1, \cP_{\varphi})}

\begin{document}

\title[Holomorphic line bundles on the loop space]
{Holomorphic line bundles on the loop space of the Riemann
sphere}
\author{Ning Zhang}
\thanks{This research was partially supported by an NSF grant.}
\address{Department of Mathematics\\ Purdue University\\ West Lafayette, IN 47906}
\email{nzhang@math.purdue.edu} \keywords{Picard group, Dolbeault
cohomology group, Loop space} \subjclass{58B12, 32Q99, 58D15}
\begin{abstract}
The loop space $L\bP_1$ of the Riemann sphere consisting of all
$C^k$ or Sobolev $W^{k,p}$ maps $S^1 \to \bP_1$ is an infinite
dimensional complex manifold. The loop group $L\pgl$ acts on
$L\bP_1$. We prove that the group of $L\pgl$ invariant holomorphic
line bundles on $L\bP_1$ is isomorphic to an infinite dimensional
Lie group. Further, we prove that the space of holomorphic
sections of these bundles is finite dimensional, and compute the
dimension for a generic bundle.
\end{abstract}
\maketitle

\section{Introduction}
Let $M$ be a finite dimensional complex manifold. Its loop space
$LM$ with a specified regularity, for example $C^k$ $(1\le k\le
\infty)$ or $W^{k,p}$ $(1\le k < \infty, 1\le p<\infty)$, consists
of all maps of the circle $S^1$ into $M$ with the given
regularity. $LM$ is an infinite dimensional complex manifold. This
paper studies holomorphic line bundles on the loop space $L\bP_1$
of the Riemann sphere.

A direct motivation comes from \cite{mz}, where Millson and Zombro
conjecture that there exists a $\pgl$ equivariant embedding of
$L\mathbb{P}_1$ into a projectivized Banach/Fr\'echet space. The
conjecture arises in connection with extending Mumford's geometric
invariant theory to an infinite dimensional setting. Another
indirect motivation comes from \cite{w}, where Witten suggests to
study the geometry and analysis of real and complex manifolds
through their loop spaces. In finite dimensions it is a problem of
fundamental importance to identify the Picard group of holomorphic
line bundles on a complex manifold and the space of holomorphic
sections of these bundles. Here we address this problem for a
class of holomorphic line bundles on the first interesting loop
space: $L\bP_1$, and in particular make some progress toward
answering the conjecture by Millson and Zombro.

The following are the main results of this paper.

If a group $G$ acts on a set $V$, let $V^G$ denote the $G$-fixed
subset of $V$. The loop space of a finite dimensional complex Lie
group is a complex Lie group under pointwise group operation (loop
group). Let Pic($L\bP_1$) be the Picard group of $L\bP_1$. The
group $\pgl$ acts on $\mathbb{P}_1$, so the loop group $L\pgl$
acts on $L\mathbb{P}_1$ and on Pic($L\bP_1$). Let
$L\mathbb{C}^{\ast}$ be the loop group of
$\mathbb{C}^{\ast}=\mathbb{C}\setminus\{0\}$, and
Hom($L\mathbb{C}^{\ast},\mathbb{C}^{\ast}$) be the group of
holomorphic homomorphisms from $ L\mathbb{C}^{\ast}$ to
$\mathbb{C}^{\ast}$.

\begin{thm}
$\Pic(L\bP_1)^{L\pgl} \cong \Hom(L\bC^{\ast}, \mathbb{C}^{\ast})$
as groups.
\end{thm}

Note that $\h$ is a $\bZ$-module of infinite rank, while the group
of topological isomorphism classes of line bundles on $L\bP_1$ is
isomorphic to $\bZ$ (cf. \cite{mc}).

Evaluation of loops in $\bP_1$ at $t\in S^1$ gives rise to a
holomorphic map $E_t: L\bP_1 \rightarrow \bP_1$. Let $\cP_t\in
\Pic(L\bP_1)$ be the pull back of the hyperplane bundle on $\bP_1$
by $E_t$.

\begin{thm}
Let $\cP \in \Pic(L\bP_1)^{L\pgl}$.
\begin{enumerate}
\item[\textup{(1)}]
If $\cP \cong \cP_{t_1}^{n_1} \otimes \cdots \otimes
\cP_{t_r}^{n_r}$, where $n_i \ge 0$ and $t_i \neq t_j$ for $i\neq
j$, then $(n_1+1)\cdots(n_r+1) \le \dim H^0(L\bP_1,\cP) <\infty$.
\item[\textup{(2)}] Otherwise $H^0(L\bP_1, \cP)=0$.
\end{enumerate}
\end{thm}

Therefore the sections of no $L\pgl$ invariant holomorphic line
bundle will give rise to a projective embedding of $L\bP_1$.

The isomorphism in Theorem 1.1 is gotten by an explicit
construction in Section 2. In Section 3 we prove Theorem 1.2(2).
In Section 4 we study the space of holomorphic sections or the
zero order Dolbeault cohomology group of line bundles defined in
Theorem 1.2(1), and in particular prove Theorem 1.2(1). Any such
line bundle obviously has holomorphic sections: products of pulled
back sections by the evaluation maps. We will show that for a
generic bundle of this type these are all sections. Yet there are
bundles which have other sections as well; interestingly, in this
case $\dim\s$ depends on the regularity of the loops.

It is natural to ask whether
$\Pic(L\bP_1)^{\pgl}=\Pic(L\bP_1)^{L\pgl}$. If so, then the
conjecture made by Millson and Zombro is answered in the negative.

The author would like to thank Professor L\'aszl\'o Lempert for
many helpful conversations and his continuous encouragement.
Thanks are also due to Professor John J. Millson and Professor
James McClure for valuable communications on some issues related
to this paper.

\section{Identification of $\hi$}

We fix a regularity class $\cF$ among $C^k$ $(1\le k\le \infty)$
respectively $W^{k,p}$ $(1\le k < \infty, 1\le p<\infty)$. In this
paper we write $LM$ ($L_kM$ resp. $L_{k,p}M$) to denote the $\cF$
($C^k$ resp. $W^{k,p}$) loop space of a manifold $M$. Let $M$ and
$N$ be finite dimensional complex manifolds and $\phi :
M\rightarrow N$ be a holomorphic map. Define $L \phi: LM \ni x
\mapsto \phi\circ x \in LN$. Then $LM$ and $LN$ are infinite
dimensional complex manifolds locally biholomorphic to open
subsets of complex Banach (Fr\'echet when $\cF=C^{\infty}$)
spaces, and $L\phi$ is holomorphic. Thus $L$ is a functor from the
category of finite dimensional complex manifolds to the category
of all complex manifolds. Let $t\in S^1$. The evaluation map
$E_t=E_t^{LM}: LM \ni x \mapsto x(t)\in M$ is holomorphic. See
Section 2 of \cite{l2}.

We call constant maps $S^1 \to M$ point loops in $M$. They form a
submanifold of $LM$, which we identify with $M$.

Next we define a map $\lb: \Hom(L\bC^{\ast}, \mathbb{C}^{\ast})
\to \hi$. We will show that $\lb$ is an isomorphism of groups,
which will then prove Theorem 1.1.

In Section 6 of \cite{mz} Millson and Zombro construct a
holomorphic line bundle on $L\bP_1$, and a similar idea in fact
yields a map from $\Hom(L\bC^{\ast}, \mathbb{C}^{\ast})$ to
$\Pic(L\bP_1)$ as follows. Let $p: Q \to \bP_1$ be the principal
$\bC^{\ast}$-bundle associated with the hyperplane bundle $H \to
\bP_1$. Applying the loop functor we obtain a principal
$L\bC^{\ast}$-bundle $Lp: LQ \to L\bP_1$. Now a homomorphism
$\varphi: L\bC^{\ast} \to \bC^{\ast}$ determines a representation
of $L\bC^{\ast}$ on $\bC$. Recall that, in general, with a
principal $G$-bundle $P\to B$ and a representation $\rho$ of $G$
on a vector space $V$, one can functorially associate a vector
bundle $E\to B$ with typical fiber $V$ (see Section 12.5 of
\cite{hu}). If $h_{ab}$ are the $G$-valued transition functions of
$P$ with respect to some trivialization, the corresponding
transition functions of $E$ will be $\rho(h_{ab})$. Accordingly we
associate with $Lp$ and $\varphi$ a line bundle $\lp$. Define the
map
$$\lb: \Hom(L\bC^{\ast},
\mathbb{C}^{\ast}) \to \Pic(L\bP_1), \hspace{1mm} \varphi \mapsto
\lp.$$

Note that the $\pgl$ action on $\bP_1$ can be covered by a $GL(2,
\bC)$ action on $Q$. One way to see this is to pass to the
tautological $\bC^{\ast}$-bundle $Q^{-1}$, whose total space is
$\bC^2\setminus \{0\}$, on which the $GL(2,\bC)$ action is
standard. The $GL(2,\bC)$ action on $Q$ gives rise to an
$LGL(2,\bC)$ action on $LQ$. Since $LGL(2,\bC) \to L\pgl$ is
surjective (as follows from the exact homotopy sequence associated
with the fibration $\bC^{\ast} \to GL(2,\bC) \to \pgl$ and the
lifting of homotopies), this $LGL(2,\bC)$ action will in fact
cover the $L\pgl$ action on $L\bP_1$. In particular,
$\gamma^{\ast}LQ \cong LQ$ for $\gamma \in L\pgl$. Hence
$\gamma^{\ast}\lp \cong \lp$ and we have proved the following

\begin{prop}
The range of $\lb$ is in $\Pic(L\bP_1)^{L\pgl}$.
\end{prop}

Let $\fU=\{U_a=\bP_1\setminus \{a\}: a\in \bP_1\}$. Then
$$L\fU=\{LU_a: a\in \bP_1\} \leqno{(2.1)}$$ is an open covering of $L\bP_1$. Now
we introduce a way to construct \v{C}ech cohomology classes in
$H^1(L\fU, \co^{\ast})$. Let $\co^G$ denote the sheaf of
holomorphic maps to the complex Lie group $G$ from a complex
manifold. So $\co^{\bC^{\ast}}=\co^{\ast}$. If $c=(c_{ab})$ is a
\v{C}ech 1-cocycle of $\fU$ with values in the sheaf $\co^{\ast}$
and $[c]$ its cohomology class, then $Lc=(Lc_{ab})$ is a \v{C}ech
1-cocycle of $L\fU$ with values in $\co^{L\bC^{\ast}}$. Any
$\varphi \in \h$ induces a sheaf homomorphism $\co^{L\bC^{\ast}}
\to \co^{\ast}$. Since the cohomology class of $\varphi \circ Lc$
depends only on $[c]$, we obtain a map
$$H^1(\fU, \co^{\ast}) \times \h \to H^1(L\fU, \co^{\ast}),
\hspace{1mm} ([c], \varphi) \mapsto [\varphi \circ Lc],
\leqno{(2.2)}$$ a group homomorphism in both variables.

Fix $[c]$ in (2.2) to be the class $[c_H]$ of the hyperplane
bundle $H \to \bP_1$, where
\begin{align*}
& c_H=\{g_{ab} \in \co^{\ast}(U_a\cap U_b): a,b \in
\bP_1, \hspace{1mm} a\neq b \}, \\
&g_{ab}=\left\{
\begin{array}{r@{,\quad}l}
\frac{z-b}{z-a}&a,b\not=\infty\\
z-b&a=\infty\\
\frac{1}{z-a}&b=\infty \end{array} \right., \hspace{1mm}z\in U_a
\cap U_b. \tag{2.3}
\end{align*}
Then we obtain a homomorphism $\h \to H^1(L\fU, \co^{\ast})$,
$\varphi \mapsto [c_{\varphi}]$, where
$$c_{\varphi}=\{\varphi\circ Lg_{ab}\in
\co^{\ast}(LU_a \cap LU_b): \hspace{1mm} a,b \in \bP_1,
\hspace{1mm} a\neq b\}, \leqno{(2.4)}$$ which has the same range
as the map in (2.2).

\begin{prop} \label{iso}
The line bundle associated to $[c_{\varphi}]$ is $\lp$. In
particular $\lb$ is a homomorphism.
\end{prop}

\begin{proof}
From the definitions it follows that $Lc_H=(Lg_{ab})$ is a family
of transition functions of $Lp$, hence $c_{\varphi}$ in (2.4) is a
family of transition functions of $\lp$.
\end{proof}

Let $t\in S^1$ and $\cP_t \in \Pic(L\bP_1)$ be the pull back of
the hyperplane bundle on $\bP_1$ by the evaluation map
$E_t^{L\bP_1}$. Obviously $E_t^{L\bC^{\ast}} \in \Hom(L\bC^{\ast},
\mathbb{C}^{\ast})$, and by Proposition \ref{iso}
$\lb(E_t^{L\bC^{\ast}})=\cP_t$.

To prove Theorem 1.1 we need the following preparations.

Identify $LU_{\infty}$ with $L\bC$ and let $y\in L\bC^{\ast}$.
Define the map $\Phi_y: L\bC\times \bP_1\ni (x, \lambda) \mapsto
x+\lambda y \in L\bP_1$, where by $\Phi_y(x,\infty)$ we mean the
point loop $\infty \in L\bP_1$. Since a continuous map $h:
\Omega_1 \to \Omega_2$ between open subsets of Fr\'echet spaces is
holomorphic if and only if its restriction to the intersection of
$\Omega_1$ with any affine line is holomorphic, see Sections 2.3,
3.1 of \cite{he}, one can easily check that $\Phi_y$ is
holomorphic. Clearly $\Phi_y(x,\cdot): \bP_1\to L\bP_1$ is a
holomorphic embedding, whose image we will denote by $v(x,y)$. The
existence of subvarieties $v(x,y)$ immediately implies that any
holomorphic function on $L\bP_1$ is constant.

Let $\cP \in \Pic(L\bP_1)$. If $\cP|_{\bP_1} \cong H^n$, then we
say that $\cP$ is of order $n$, or $\ord(\cP)=n$; and we claim
that $\cP|_{v(x,y)} \cong H^n$, $x \in L\bC$, $y \in L\bC^{\ast}$.
This would imply that the only holomorphic section of $\cP$ is the
zero section if $\ord(\cP)< 0$. To show the claim, let $c_1$
denote the rational first Chern class of a line bundle. According
to the first theorem of \cite{vps}, the inclusion $\bP_1 \to
L\bP_1$ induces an isomorphism $H^2(L\bP_1, \bQ) \cong H^2(\bP_1,
\bQ)$, so that $c_1(\cP)$ is completely determined by
$c_1(\cP|_{\bP_1})=n$. Hence $c_1(\cP|_{v(x,y)})$ , and therefore
the degree of $\cP|_{v(x,y)}$, is also determined by $n$. That
this degree is itself $n$ then follows from computing it in the
special case $\cP=\cP^n_t$.

Let $\varphi \in \Hom(L\bC^{\ast}, \bC^{\ast})$. The restriction
of $\varphi$ to the subgroup of point loops must be of the type $z
\mapsto z^n$, $z\in \bC^{\ast}$, where $n\in \bZ$. We call this
$n$ the order of $\varphi$ and denote it by $\op$. Proposition
\ref{iso} implies that $\ord(\lp)=\op$. Let $F_{x}$ be the fiber
of $\cP$ at $x\in L\bP_1$.

\begin{prop} \label{tri}
Let $\cP \in \Pic(L\bP_1)$, $\ord(\cP)=0$ and $y \in L\bC^{\ast}$.
For any $\zeta \in F_{\infty}\! \setminus \! \{0\}$ there exists a
unique non-vanishing section $\sigma=\sigma_{y, \zeta} \in
H^0(LU_{\infty}, \cP)$ such that $\lim_{\lambda \to
\infty}\sigma(x+\lambda y)=\zeta$, $x \in LU_{\infty}$. In
particular, $\cP|_{LU_{\infty}}$ is holomorphically trivial.
\end{prop}

\begin{proof}
Since $\cP|_{v(x,y)}$ is trivial, uniqueness is obvious. As to
existence, let $\cH$ be a hyperplane in $L\bC$ such that $y \notin
\cH$ and consider the line bundle
$\tilde{\cP}=\Phi_y^{\ast}\cP|_{\cH \times \bP_1}$. Since $\Phi_y
\equiv \infty$ on $\cH \times \{\infty\}$, $s=\Phi_y^{\ast}\zeta$
is a non-vanishing holomorphic section of $\tilde{\cP}|_{\cH
\times \{\infty\}}$. In turn $s$ determines a section
$\tilde{\sigma}$ of $\tilde{\cP}|_{\cH \times \bP_1}$ such that
$\tilde{\sigma}|_{\{x\} \times \bP_1}$ is the unique non-vanishing
holomorphic section of $\tilde{\cP}|_{\{x\} \times \bP_1}$ with
$\tilde{\sigma}(x,\infty)=s(x,\infty)$, for all $x \in \cH$. Next
we show that $\tilde{\sigma}$ is holomorphic. Let $x_0 \in \cH$.
By Proposition 5.1 of \cite{l1} there exist a neighborhood $x_0
\in U \subset \cH$ and a section $v \in C^{\infty}(U\times \bP_1,
\tilde{\cP})$ such that $v$ is holomorphic on $\{x\}\times \bP_1$
for all $x \in U$, and $v|_{\{x_0\} \times
\bP_1}=\tilde{\sigma}|_{\{x_0\} \times \bP_1} \neq 0$. By choosing
a sufficiently small $U$ we can assume $v \neq 0$. The function
$\tilde{\sigma}/v$ is $C^{\infty}$ on $U \times \{\infty\}$ and is
constant on $\{x\}\times \bP_1$, $x \in U$, hence is $C^{\infty}$
on $U\times \bP_1$. So it follows that $\tilde{\sigma} \in
C^{\infty}(\cH \times \bP_1)$. Since $\ba \tilde{\sigma}|_{\cH
\times \{\infty\}}=0$ and $\ba \tilde{\sigma}|_{\{x\} \times
\bP_1}=0$ for all $x \in \cH$, by Proposition 5.2(ii) of \cite{l1}
we obtain that indeed $\ba \tilde{\sigma}=0$. Then the desired
$\sigma$ is the pull back of $\tilde{\sigma}$ by
$(\Phi_y|_{\cH\times \bC})^{-1}$.
\end{proof}

Since $H|_{U_{\infty}}$ is trivial, so is $\cP_t|_{LU_{\infty}}$.
In general, let $\cP \in \Pic(L\bP_1)$, $\ord(\cP)=n$. As
$$\cP=\cP_t^n \otimes (\cP^{-n}_t \otimes \cP), \hspace{1mm}
\mathrm{where} \hspace{1.5mm} \ord(\cP^{-n}_t \otimes \cP)=0,
\leqno{(2.5)}$$ Proposition \ref{tri} implies that
$\cP|_{LU_{\infty}}$ is also trivial. More generally,
$\cP|_{LU_a}$ is trivial, $a \in \bP_1$, which means

\begin{cor}
$\Pic(L\bP_1) \cong H^1(L\fU, \co^{\ast})$.
\end{cor}

If $\cP\in\hi$, then Proposition \ref{tri} can be improved:
$\sigma$ there is essentially independent of $y$, and so is a
canonical section of $\cP|_{LU_{\infty}}$.

\begin{prop} \label{can}
Let $\cP \in \hi$, $\ord(\cP)\ge 0$. Then there exists a
non-vanishing section $\sigma_{\infty} \in H^0(LU_{\infty}, \cP)$
such that $\lim_{\lambda \to \infty}\sigma_{\infty}(x+\lambda y)$
exists for all $x \in LU_{\infty}$ and $y\in L\bC^{\ast}$. Such a
section is unique up to a multiplicative constant.
\end{prop}

\begin{proof}
Suppose $\sigma_{\infty}$ exists. For any $y \in L\bC^{\ast}$,
$\cP|_{v(0,y)}\cong H^{\ord(\cP)}$ has a unique holomorphic
section which does not vanish on $v(0,y) \! \setminus \!
\{\infty\}$ and assumes $\sigma_{\infty}(0)$ at $0$; clearly
$\sigma_{\infty}$ agrees with this section on $v(0,y) \! \setminus
\! \{\infty\}$. In particular, $\sigma_{\infty}(y)$ is uniquely
determined by $\sigma_{\infty}(0)$, $y\in L\bC^{\ast}$. Since
$L\bC^{\ast}$ is dense in $L\bC=LU_{\infty}$, $\sigma_{\infty}$ is
completely determined by its value at $0$, or unique up to a
multiplicative constant. Next we show the existence.

First assume $\ord(\cP)=0$. Let $\zeta \in F_{\infty}\! \setminus
\! \{0\}$, $y_1, y_2 \in L\bC^{\ast}$, and $\sigma_{y_1, \zeta},
\sigma_{y_2, \zeta} \in H^0(LU_{\infty}, \cP)$ as in Proposition
\ref{tri}. Define $$h=h_{y_1, y_2}=\frac{\sigma_{y_1,
\zeta}}{\sigma_{y_2, \zeta}} \in \co^{\ast}(LU_{\infty}),$$ which
is independent of the choice of $\zeta \in F_{\infty} \! \setminus
\! \{0\}$. For fixed $x_1 \in L\bC$ let $\mu \in L\pgl$ be the
translation $x \mapsto x+x_1$, $x \in L\bP_1$. Clearly
$\mu^{\ast}\sigma_{y_i, \zeta} \in H^0(LU_{\infty},
\mu^{\ast}\cP)$, $i=1,2$, is a section of the type as in
Proposition \ref{tri}. With a fixed isomorphism $\mu^{\ast}\cP
\cong \cP$, $\mu^{\ast}\sigma_{y_i, \zeta}$ corresponds to
$\sigma_{y_i, \zeta^{\prime}}$, where $\zeta^{\prime} \in
F_{\infty} \! \setminus \! \{0\}$. Therefore
$$h(x+x_1)=\mu^{\ast}h(x)=\frac{\mu^{\ast}\sigma_{y_1,
\zeta}}{\mu^{\ast}\sigma_{y_2, \zeta}}(x)=\frac{\sigma_{y_1,
\zeta^{\prime}}}{\sigma_{y_2, \zeta^{\prime}}}(x)=h(x).$$ Thus
$h=h_{y_1, y_2}$ is a non-zero constant. From the definition of
$\sigma_{y_2, \zeta}$ it follows that $\sigma_{y_1,
\zeta}=h\sigma_{y_2, \zeta}$ extends to be a section of
$\Lambda|_{v(x, y_2)}$ for all $x \in L\bC$, where $y_2\in
L\bC^{\ast}$ is arbitrary. The upshot is
$\sigma_{\infty}=\sigma_{y_1, \zeta}$ with arbitrary $y_1, \zeta$
will do.

Second assume $\cP=\cP_t^n$, $n\ge 0$. If $s$ is a section of
$H^n$, nonzero on $\bP_1\setminus \{\infty\}$, take
$\sigma_{\infty}=(E_t^{L\bP_1})^{\ast}s$. Finally, these two
special cases and (2.5) imply the general case.
\end{proof}

Theorem 1.1 follows from

\begin{thm} \label{pi}
$\lb: \Hom(L\bC^{\ast}, \mathbb{C}^{\ast}) \to
\Pic(L\bP_1)^{L\pgl}$ is a group isomorphism.
\end{thm}

\begin{proof}
i) Injectivity. If $\lb(\varphi)=\lp$ is trivial, then $\op=0$;
and by Proposition \ref{iso} and (2.4) we can find $f_a \in
\co^{\ast}(LU_a)$, $a \in \bP_1$, such that $f_a=(\varphi \circ
Lg_{ab})f_b$ on $LU_a\cap LU_b$. In particular,
$f_{\infty}=\varphi f_0$ on $ L\bC^{\ast} \subset L\bP_1$. For any
$y \in L\bC^{\ast}$ we have
$$\lim_{\bC \ni \lambda \to \infty}f_{\infty}(\lambda
y)=\lim_{\lambda \to \infty}\varphi(\lambda y)f_0(\lambda
y)=\varphi(y)f_0(\infty),$$ since $\varphi(\lambda
y)=\varphi(\lambda)\varphi(y)=\varphi(y)$. Thus
$f_{\infty}|_{v(0,y)\setminus \{\infty\}}$ extends to all of
$v(0,y)$, and so must be constant. In particular,
$f_{\infty}(y)=f_{\infty}(0)$ and $f_{\infty}$ itself is a
constant. Similarly $f_0$ is also a constant. Then $\varphi$ is a
constant which can only be $1$.

ii) Surjectivity. Since $\cP_t^n$ is in the range of $\lb$, by
(2.5) we only need to show that any $\cP \in \hi$ with
$\ord(\cP)=0$ is in the range of $\lb$.

Let $\varepsilon(x)=1/x$, $x \in L\bP_1$. The induced bundle
$\varepsilon^{\ast}\cP$ is isomorphic to $\cP$, so Proposition 2.5
applies to produce a non-vanishing $\tilde{\sigma}_{\infty} \in
H^0(LU_{\infty}, \varepsilon^{\ast}\cP)$. Then
$\sigma_0=\tilde{\sigma}_{\infty} \circ \varepsilon \in H^0(LU_0,
\cP)$ is characterized, up to a multiplicative constant, by the
fact that $\lim_{\lambda \to \infty}\sigma_0\left((x+\lambda
y)^{-1}\right)$ exists, for all $x \in LU_{\infty}$, $y \in
L\bC^{\ast}$. To get rid of the ambiguity in the choice of the
constant, fix a non-zero $s \in H^0(\bP_1, \cP)$ and choose
$\sigma_0$ and $\sigma_{\infty}$ (as in Proposition 2.5) to agree
with $s$ on point loops. Set
$\phi=\sigma_0/\sigma_{\infty}:L\bC^{\ast} \to \bC^{\ast}$. Note
that on point loops $\phi=1$. We will show that $\phi$ is a
homomorphism and $\cP=\cP_{\phi}$.

For this purpose fix $y_1 \in L\bC^{\ast}$ and define
$\gamma(x)=y_1x$, $x \in L\bP_1$. It is straightforward that the
non-vanishing section $\gamma^{\ast}\sigma_{\infty} \in
H^0(LU_{\infty}, \gamma^{\ast}\cP)$ satisfies the conditions in
Proposition 2.5. Hence under an isomorphism $\gamma^{\ast}\cP
\cong \cP$, $\gamma^{\ast}\sigma_{\infty}$ corresponds to a
constant multiple of $\sigma_{\infty}$. Similarly, under this
isomorphism $\gamma^{\ast}\sigma_0$ corresponds to a constant
multiple of $\sigma_0$. Therefore
$$\phi(y_1 y)=\left(\gamma^
{\ast}\phi\right)(y)
=\frac{\gamma^{\ast}\sigma_0}{\gamma^{\ast}\sigma_{\infty}}(y)=
c\frac{\sigma_0}{\sigma_{\infty}}(y)=c\phi(y).$$ Letting $y=1$ we
get $c=\phi(y_1)$, so $\phi$ is indeed a homomorphism.

Finally with $a \in \bC$ let $\mu_a(x)=x+a, x \in L\bP_1$. For the
bundle $\cP^{\prime}=\mu_a^{\ast}\cP$ one can construct
corresponding sections $\sigma^{\prime}_{\infty}$ and
$\sigma^{\prime}_0$; for the normalization, use
$s^{\prime}=\mu_a^{\ast} s \in H^0(\bP_1, \cP^{\prime})$. As
$\cP^{\prime} \cong \cP$,
$\sigma_0^{\prime}/\sigma_{\infty}^{\prime}=\phi$. Since
$(\mu_a^{-1})^{\ast}\sigma_{\infty}^{\prime} \in H^0(LU_{\infty},
\cP)$ satisfies the conditions in Proposition 2.5 and agrees with
$s$ on point loops,
$(\mu_a^{-1})^{\ast}\sigma_{\infty}^{\prime}=\sigma_{\infty}$. Let
$\sigma_a=(\mu_a^{-1})^{\ast}\sigma_{0}^{\prime} \in H^0(LU_a,
\cP)$. Since
$$\frac{\sigma_a}{\sigma_{\infty}}=\left(\mu_a^{-1}\right)^{\ast}
\frac{\sigma_0^{\prime}}{\sigma_{\infty}^{\prime}}=\phi\circ
\mu_a^{-1},$$ it is straightforward to check that the transition
functions $\sigma_b/\sigma_a \in \co^{\ast}(LU_a \cap LU_b)$ of
$\cP$ agree with those of $\cP_{\phi}$ given in (2.4); therefore
$\cP=\cP_{\phi}$ is indeed in the range of $\lb$.
\end{proof}

Let $(L\bC)^{\ast}$ be the space of continuous linear functionals
on $L\bC$, let $\tp \in (L\bC)^{\ast}$ be the Lie algebra
homomorphism induced by $\varphi \in \h$, and $x_0 \in
L\bC^{\ast}$ be a fixed loop whose winding number with respect to
$0$ is $1$. The reader can check that the map
$$\h \to \{\phi \in (L\bC)^{\ast}: \phi(1)\in\bZ \}\times \bC^{\ast},
\hspace{1mm} \varphi \mapsto (\tp, \varphi(x_0))$$ is an
isomorphism of groups. Therefore $\h$ is a $\bZ$-module of
infinite rank.

The proof of Theorem 2.5 implies the following

\begin{prop}
Let $\varphi \in \h$, $\ord(\varphi) \ge 0$. There is a family of
non-vanishing sections $\{\sigma_a \in H^0(LU_a, \lp): a\in
\bP_1\}$, unique up to an overall multiplicative constant, that
satisfies $\sigma_b/\sigma_a=\varphi \circ Lg_{ab}$ on $LU_a \cap
LU_b$. Furthermore, $\sigma_{\infty}$ satisfies the conditions in
Proposition 2.5.
\end{prop}

\begin{proof} \label{secs}
Uniqueness is obvious: if $\{\sigma_a^{\prime}\}$ is another
family then $\{\sigma_a^{\prime}/\sigma_a\}$ defines a holomorphic
function on $L\bP_1$, which, as we have said, must be a constant.
When $\ord(\cP)=0$, the family $\{\sigma_a\}$ is constructed in
the proof of Theorem 2.6. When $\cP=\cP^n_t$, take sections
$\tau_a$ of the hyperplane bundle $H \to \bP_1$ such that
$\tau_b/\tau_a=g_{ab}$ (see (2.3)), then
$\sigma_a=(E_t^{L\bP_1})^{\ast}\tau_a^n|_{LU_a}$ will do. The case
of a general $\cP$ now follows from (2.5).
\end{proof}

\section{Proof of Theorem 1.2(2)}

We start with two results concerning polynomials on $L\bC$. For
simplicity let $E_t$ denote both $E_t^{L\bC} \in (L\bC)^{\ast}$
and $E_t^{L\bC^{\ast}} \in \mbox{Hom}(L\mathbb{C}^{\ast},
\mathbb{C}^{\ast})$ till the end of this paper.

\begin{lemma} \label{lema} If a homogeneous polynomial $h \in \co(L\bC)$ of
degree $n \ge 1$ does not vanish on $L\bC^{\ast}$, then
$h=cE_{t_1}\cdots E_{t_n}$, where $t_1, \cdots, t_n \in S^1$, and
$c \neq 0$ is a constant.
\end{lemma}
\begin{proof}
Since $L_{\infty}\bC$ is embedded into $L\bC$ with a dense image,
we only need to show the lemma for the case of $C^{\infty}$ loops.

Let $W=L\bC \setminus L\bC^{\ast}$, $\cZ_h \subset W$ be the zero
locus of $h$, and $\Emb(S^1, \mathbb{C})\subset
L_{\infty}\mathbb{C}$ be the open subset of embedded loops. Then
$\cZ_h \cap \mbox{Emb}(S^1, \mathbb{C}) \neq \emptyset$. Otherwise
for any $x \in \mbox{Emb}(S^1, \mathbb{C})$, the polynomial
$h(x+\lambda)$ in $\lambda \in \bC$ has no zero, hence is
constant. In particular, $h(x+1)=h(x)$ on $\mbox{Emb}(S^1,
\mathbb{C})$ and therefore on $L_{\infty}\bC$. So $h(1)=h(0)=0$.
But $1 \in L\bC^{\ast}$, contradiction.

Let $x_1 \in W \cap \mbox{Emb}(S^1, \mathbb{C})$. Next we show
that $W$ is a submanifold of real codimension one near $x_1$. Let
$t_1$ be the unique element of $S^1=\bR/\bZ$ such that
$x_1(t_1)=0$. We can assume that $x_1^{\prime}(t_1)$ is not real,
otherwise replace $x_1$ by $ix_1$. Consider the equation
$$\text{Im}\,x(s)=0, \hspace{1.5mm} x\in L_{\infty}\mathbb{C},\, s\in
S^1. \leqno{(3.1)}$$ This equation can also be considered on the
$C^k$ loop space $L_k\bC$, $1 \le k <\infty$. Note that
$\text{Im}\,x_1^{\prime}(t_1) \not=0$. Apply the Implicit Function
Theorem on Banach spaces (see Theorem 2.5.7 of \cite{amr}) to the
$C^k$ map $L_k\bC \times S^1 \to \bR$, $(x, t) \mapsto \text{Im}\,
x(t)$ near $(x_1,t_1)$. When $k=1$ we obtain a neighborhood $U_1
\subset L_1\bC$ of $x_1$ consisting of embedded loops, a
neighborhood $V \subset S^1$ of $t_1$, and a $C^1$ map $\phi: U
\rightarrow V$ such that for any $x \in U$, $s=\phi(x)$ is the
unique solution of (3.1) in $V$. We can shrink $U_1$ and $V$ if
necessary to ensure that $\text{Im}\,x^{\prime}(t) \neq 0$ for all
$(x,t) \in U_1 \times V$. For arbitrary $k < \infty$ we obtain
that $\phi$ is $C^k$ on $U_k=U_1 \cap L_k\bC$ by the Implicit
Function Theorem at any $(x, \phi(x))$, $x \in U_k$. This implies
that $\phi$ is $C^{\infty}$ on $U=U_1 \cap L_{\infty}\bC$. We can
choose sufficiently small $U$ so that $x(t)\neq 0$ if $x\in U$ and
$t\not\in V$; then
$$x\left(\phi(x)\right)=0, \hspace{1mm} x \in U\cap W. \leqno{(3.2)}$$
Let $Y$ be the real hyperplane $\{y \in L_{\infty}\mathbb{C}:
\text{Im}\,y(t_1)=0\}$ of $L_{\infty}\bC$. Define the $C^{\infty}$
map $\tau: S^1 \times Y \rightarrow L_{\infty}\mathbb{C}$, $(s,
x(t)) \mapsto x(t-s)$. The $C^{\infty}$ map $\rho: U \rightarrow
S^1 \times Y$, $x(t)\mapsto (\phi(x)-t_1, x(t+\phi(x)-t_1))$ is
the local inverse of $\tau$ near $(t_1, x_1)$. Let $\cH_{t}$ be
the kernel of $E_{t}$. Note that $\cH_{t_1} \subset Y$. Since
$\tau(S^1 \times \cH_{t_1})=W$ and by (3.2) $\rho(W \cap U)
\subset S^1 \times \cH_{t_1}$, $W \cap U$ is a submanifold of real
codimension one of $U$. Its tangent space is
$$TW_{x}=\cH_{\phi(x)} \oplus \{x^{\prime}\mathbb{R}\}, \hspace{1mm} x \in W \cap
U,\leqno{(3.3)}$$ if we identify the tangent space of
$L_{\infty}\bC$ at $x$ with $L_{\infty}\mathbb{C}$.

Now assume $x_1 \in \mbox{Emb}(S^1, \mathbb{C}) \cap
\mathcal{Z}_h$. Let $v \in \cH_{t_1}$ and $v_0 \not\in W$.
Consider the restriction of $h$ to the $2$-dimensional affine
subspace $E=\{x_1+\lambda_1 v+\lambda_2 v_0: \lambda_1, \lambda_2
\in \mathbb{C}\}$. Let $\delta: \Delta \to \cZ_h \cap E \cap U
\subset W$ be a holomorphic arc such that $\delta(0)=x_1$ and $0
\neq \delta^{\prime}(0) \in TW_{x_1}$. As for all $\lambda \in
\Delta$ $\delta^{\prime}(\lambda)$ is in the maximal complex
subspace of $T_{\delta(\lambda)}W$, $\delta^{\prime}(\lambda) \in
\cH_{\phi(\delta(\lambda))}$ by (3.3). Since $\phi$ is constant on
$\cH_{\phi(\delta(\lambda))} \cap U$, the derivative of $\phi$ in
the direction of $\delta^{\prime}(\lambda)$ is zero, hence $\phi
\circ \delta \equiv t_1$. In view of (3.2) $\delta(\Delta) \subset
\cH_{t_1} \cap E$. It follows that $\{x_1+\lambda_1 v\} \subset
\cZ_h$ and $\cH_{t_1} \subset \mathcal{Z}_h$.

The local ring $\mathcal{O}(L_{\infty}\mathbb{C})_{x}$, $x\in
L_{\infty}\bC$, is a unique factorization domain, see Proposition
5.15 of \cite{ma}. The germ of $E_{t_1}$ at any $x \in \cH_{t_1}$
is prime, for the functions $e(\lambda)=E_{t_1}(x+\lambda y)$
vanish to first order at $\lambda =0$ if $y(t_1) \neq 0$. Applying
the Nullstellensatz in Theorem 5.14 of \cite{ma} to the prime
ideals $((E_{t_1})_{x})$, $x \in \cH_{t_1}$, we obtain that
$h=E_{t_1}\tilde{h}$, $\tilde{h} \in
\mathcal{O}(L_{\infty}\mathbb{C})$. On any affine line $\tilde{h}$
is a polynomial of degree $\le n-1$, so it is a polynomial of
degree $\le n-1$ on $L_{\infty}\mathbb{C}$, see Section 2.2 of
\cite{he}. The zero locus of $\tilde{h}$ is still in $W$, so that
repeating the above process we obtain the conclusion of the lemma.
\end{proof}

Let $\cL^n(L\bC)$ be the space of holomorphic polynomials of
degree $\le n$ on $L\bC$.

\begin{prop} \label{pair} Let $h_1, h_2 \in \co(L\bC)$ and
$\varphi \in \co(L\bC^{\ast})$ satisfy $\varphi(\lambda
y)=\lambda^n \varphi(y)$, $\lambda \in \bC^{\ast}$, $y\in
L\bC^{\ast}$, $n \in \bN \cup \{0\}$. If
$h_1(y)=\varphi(y)h_2(y^{-1})$, $y\in L\bC^{\ast}$, then $h_1,h_2
\in \cL^n(L\bC)$. Let $h_j^i$ be the $i$-th order homogeneous
component of $h_j$, $j=1,2$, $i=0, \cdots, n$. Then
$h_1^i(y)=\varphi(y)h_2^{n-i}(y^{-1})$.
\end{prop}

\begin{proof}
If $x \in L\bC$, $y \in L\bC^{\ast}$, then $$ \lim_{\lambda \to
\infty} \frac{h_1(x+\lambda y)}{\lambda^n} = \lim_{\lambda \to
\infty} \varphi(x \lambda^{-1}+y) h_2((x+\lambda y)^{-1}) =
\varphi(y)h_2(0). $$ Thus $h_1(x+\lambda y)$ is a polynomial of
degree $\le n$ in $\lambda$. As $L\bC^{\ast}$ is dense in $L\bC$,
the same holds for all $y \in L\bC$. Since a continuous function
on a Fr\'echet space is a polynomial of degree $\le n$ if its
restriction to any affine line is such a polynomial, see Section
2.2 of \cite{he}, we conclude that $h_1 \in \cL^n(L\bC)$.
Similarly $h_2\in \cL^n(L\bC)$. Comparing homogeneous components
of same order on both sides of the equation
$h_1(y)=\varphi(y)h_2(y^{-1})$, we get
$h_1^i(y)=\varphi(y)h_2^{n-i}(y^{-1})$, $i=0, \cdots, n$.
\end{proof}

Let $\sigma \in H^0(L\bP_1, \Lambda_{\varphi})$, $\ord(\varphi)
\ge 0$. With sections $\sigma_a \in H^0(LU_a, \lp), a\in \bP_1$,
in Proposition 2.7, the functions $H_a= \sigma/\sigma_a \in
\mathcal{O}(LU_a)$ satisfy
$$H_a(x)=\varphi \circ Lg_{ab}(x)H_b(x), \hspace{1mm} x\in LU_a \cap LU_b.
\leqno{(3.4)} $$ Note that $Lg_{ab}$ maps $LU_a$ to $L\mathbb{C}$
biholomorphically, and let $h_{ab}=H_a\circ Lg_{ab}^{-1} \in
\mathcal{O}(L\bC)$. Then (3.4) implies that
$$h_{ab}(y)=\varphi(y) h_{ba}(y^{-1}), \hspace{1mm}  y\in L\bC^{\ast}.
\leqno{(3.5)}$$ \noindent From Proposition \ref{pair} it follows
that $h_{ab} \in \cL^n(L\bC)$; also $$H_{\infty}=h_{\infty 0} \in
\cL^n(L\bC). \leqno{(3.6)}$$

\medskip \noindent {\it Proof of Theorem 1.2(2).} By Theorem \ref{pi}
each element of $\Pic(L\bP_1)^{L\pgl}$ is of the type $\lp$,
$\varphi \in \h$. Suppose $\lp$ has a non-zero holomorphic section
$\sigma$. We can assume that $\sigma(\infty) \neq 0$; this can be
arranged by pulling $\sigma$ back by a suitable element of
$L\pgl$. Then $\lp|_{\bP_1} \cong H^n$ also has a non-zero
section, so $\ord(\varphi)=n\ge 0$. Note that $h_{0
\infty}(0)=H_0(\infty) \neq 0$. By (3.5), applied with $ab=\infty
0$, and Proposition \ref{pair}
$$h^n_{\infty 0}(y)=\varphi(y)h_{0 \infty}(0) \neq 0, \hspace{1mm} y \in
L\bC^{\ast}. \leqno{(3.7)}$$ Hence $h=h_{\infty 0}^n \in
\co(L\bC)$ satisfies the condition in Lemma \ref{lema}, so
$h^n_{\infty 0}=cE_{t_1}\cdots E_{t_n}$. By (3.7)
$\varphi=E_{t_1}\cdots E_{t_n}$, and $\lp$ is of the type as in
Theorem 1.2(1). \qed

\section{The space of holomorphic sections}

In this section we shall study the space $H^0(L\bP_1, \cP)$, where
$\cP$ is as in Theorem 1.2(1), and we shall prove Theorem 1.2(1).
Since $\dim H^0(L\bP_1, \cP)=1$ if $\cP$ is trivial, we fix
$\cP=\lp$ nontrivial, where $$\varphi=E^{n_1}_{t_1} \cdots
E^{n_r}_{t_r}, \hspace{1.5mm} n_i>0, \hspace{1mm} t_i\neq t_j
\hspace{1mm} \mathrm{if} \hspace{1mm} i\neq j. \leqno{(4.1)}$$

With $\sigma_a \in H^0(LU_a, \lp)$ as in Proposition 2.7 and
$\sigma \in H^0(L\bP_1, \lp)$, we have seen that
$\sigma/\sigma_{\infty}=H_{\infty} \in \cL^n(L\bC)$, cf. (3.6).
Define a monomorphism
$$\hp: H^0(L\bP_1, \lp) \to \cL^n(L\bC), \hspace{1mm} \sigma \mapsto
\sigma/\sigma_{\infty},$$ where $n=\ord(\varphi) \ge 0$. Let $\R
\subset \cL^n(L\bC)$ be the range of $\hp$. We shall study $\s$
through $\R$.

The bundle $\lp$, where $\varphi$ is as in (4.1), has non-trivial
holomorphic sections: products of pull back sections by evaluation
maps on $L\bP_1$.

\begin{prop} \label{pb}
The linearly independent functions $E^{m_1}_{t_1}\cdots
E^{m_r}_{t_r}$ are in $\R$, $0 \le m_i \le n_i, 1\le i \le r$. In
particular, $(n_1+1)\cdots (n_r+1) \le \dim \s$.
\end{prop}
\begin{proof}
Choose a basis $\{\tilde{\tau}_j: j=0, \cdots, n_1\}$ of
$H^0(\bP_1, H^{n_1})$, where $\tilde{\tau}_0$ has a zero of order
$n_1$ at $\infty \in \bP_1$ and
$\tilde{\tau}_j=z^j\tilde{\tau_0}$, $z \in \bP_1$. Let
$\tau_j=(E_{t_1}^{L\bP_1})^{\ast}\tilde{\tau_j} \in H^0(L\bP_1,
\cP^{n_1}_{t_1})$. The section $\tau_0|_{LU_{\infty}}$ satisfies
the conditions in Proposition 2.5. Therefore
$$\hp(\tau_j)=\tau_j/\tau_0=E_{t_1}^{\ast}z^j
=E_{t_1}^{j}.$$ Similarly we have pull back sections of
$\Lambda_{t_i}^{n_i}$, $i=1, \cdots, r$. By taking products of
such sections we obtain sections of $\lp$, hence
$E^{m_1}_{t_1}\cdots E^{m_r}_{t_r} \in \R, 0 \le m_i \le n_i, 1\le
i \le r$. Another way of obtaining these functions is to pull back
monomials on $\bC^r$ by the surjective map $(E_{t_1}, \dots,
E_{t_r}): L\bC \to \bC^r$, which interpretation proves the claim
of linear independence.
\end{proof}

The elements of $\R$ identified in Proposition \ref{pb} are in the
subalgebra of $\cL^n(L\bC)$ generated by evaluation functions. In
other words they are polynomials in finitely many linear
functionals. The next proposition shows that this is true in
general. For convenience, in Propositions 4.2 and 4.3 we shall
restrict our discussion to $C^k$ ($1\le k \le\infty$) loop spaces
$L_k\bP_1$. Let $x^{(\nu)}(t)$ denote the $\nu$-th derivative of
$x \in L_k\bC$ at $t \in S^1=\bR/\bZ$, $\nu \le k$. Note that the
function $x \to x^{(\nu)}(t)$ is in $(L_k\bC)^{\ast}$.

\begin{prop} \label{poly}
If $\sigma \in H^0(L_k\bP_1, \lp)$ and $P=\hp(\sigma)$, then
$P(x)$ is a polynomial in finitely many derivatives
$x^{(\nu)}(t_i)$, $0 \le \nu \le k$, $1 \le i \le r$.
\end{prop}

\begin{proof}
Let $A=\{t_1, \cdots, t_r \} \subset S^1$, $x_0 \in
L_kU_{\infty}$, and denote the $k$-jet of $x \in L_k\bP_1$ by
$j^kx$. Define
$$Z=Z(k, A, x_0)=\{x \in L_k\bP_1: j^kx|_A=j^kx_0|_A\}.$$
This is a connected complex submanifold of $L_k\bP_1$ and any
holomorphic function on it is a constant, see Sections 3, 4 of
\cite{l2}. Consider the sections $\sigma_a$ of Proposition 2.7. In
the proof of that proposition we have shown that, when
$\cP=\cP_t^n$, $\sigma_a$ can be taken to be the pullback by
$E_t^{L\bP_1}$ of sections $\tau_a$ of the hyperplane bundle $H
\to \bP_1$, $\tau_a \neq 0$ on $U_a$. When $\cP=\otimes
\cP_{t_i}^{n_i}$, $\sigma_a$ can be taken as the product of such
sections. Then $\sigma_a(x)=0$ only if $x(t_i)=a$ for some $i$; in
particular $\sigma_{\infty} \neq 0$ on $Z$. It follows that
$\sigma/\sigma_{\infty}|_Z \in \co(Z)$ is constant, i.e., $P$ is
constant on any affine subspace
$$\{x \in L_k\bC: j^kx|_A=j^kx_0|_A \}, \hspace{1mm} x_0 \in
L_k\bC. \leqno{(4.2)}$$

When $k<\infty$, the continuous linear functionals $x \to
x^{(\nu)}(t_i)$, $0 \le \nu \le k$, $1 \le i \le r$, give rise to
a surjective linear map $J: L_k\bC \to \bC^{(k+1)r}$, whose fibers
are the affine subspaces in (4.2). What we have shown above
implies that $P=J^{\ast}f$, where $f$ is a function on
$\bC^{(k+1)r}$. In fact $f$ is a holomorphic polynomial of degree
$\le n$, for with a linear right inverse $I$ to $J$ we have $f=P
\circ I$.

Consider the case when $k=\infty$. Let $P^j$ be the $j$-th order
homogeneous component of $P \in \cL^n(L_{\infty}\bC)$, $j=1,
\cdots, n$. These components are also constant on any affine
subspace in (4.2). By the definition of a homogeneous polynomial
we can find a continuous symmetric $j$-linear mapping $\Psi_j:
(L_{\infty}\bC)^j \to \bC$ such that $P^j(x)=\Psi_j(x, \cdots,
x)$. Applying the Schwartz Kernel Theorem (see Theorem 5.2.1 of
\cite{ho}) one can show that there exists a distribution $K_j$ in
the $j$ dimensional torus $T^j$ such that $$\Psi_j(x_1, \cdots,
x_j)=K_j(x_1(s_1)\cdots x_j(s_j)), \hspace{1mm}(s_1, \cdots, s_j)
\in T^j.$$ So $P^j(x)=K_j(x(s_1)\cdots x(s_j))$. The Polarization
Formula (see (2) in Section 2.2 of \cite{he})
$$K_j\left(x_1(s_1)\cdots x_j(s_j)\right)=\frac{1}{2^j j!} \sum_{\varepsilon_1,
\cdots, \varepsilon_j=\pm 1} \varepsilon_1 \cdots \varepsilon_j
P^j(\varepsilon_1x_1+\cdots+\varepsilon_j x_j)$$ and the fact that
$P^j$ depends only on $j^{\infty}x|_A$ imply that $K_j$ is
supported in $A^j\subset T^j $. Therefore $K_j$ is a (finite)
linear combination of partial derivatives at points in $A^j$, see
Theorem 2.3.4 of \cite{ho}. Hence $P^j$ and $P$ are polynomials in
finitely many $x^{(\nu)}(t_i)$, $x \in L_{\infty}\bC$.
\end{proof}

For each $i$ let $N_i$ be the order of the highest derivative
$x^{(\nu)}(t_i)$ that $P(x)$ depends on, see Proposition
\ref{poly}, and $m_i$ the degree of $P(x)$ as a polynomial of
$x^{(\cN_i)}(t_i)$. Our next task is to estimate $N_i$, $m_i$.

\begin{prop} \label{up}
$m_i \left(\cN_i+1 \right) \le n_i$, where $n_i$ is defined in
(4.1).
\end{prop}

\begin{proof}
Fix $i$. At first assume $P=\hp(\sigma)$ contains the monomial
$$c_1 \, x^{(\cN_i)}(t_i)^{m_i}, \leqno{(4.3)}$$
where $c_1 \neq 0$ is a constant. By (3.6) and (3.5) $P=h_{\infty
0}$ satisfies $h_{0 \infty}(y)=\varphi(y)h_{\infty 0}(y^{-1})$ so
that Proposition 3.2 implies that
$$h_{0\infty}^{n-m_i}(x)=\varphi(x)\left[P^{m_i}(x^{-1})\right],
\hspace{1mm} x\in L_{k}\bC^{\ast}, \leqno{(4.4)}$$ where
superscripts indicate homogeneous components of the given order.
Since $P^{m_i}$ is a polynomial in $x^{(\nu)}(t_j)$,
$P^{m_i}(x^{-1})$ is a sum of rational expressions, where the
denominators are monomials in $x(t_j)$, $1 \le j \le r$, and the
numerators are monomials in $x^{(\nu)}(t_j)$, $\nu \ge 1$. In this
sum the monomial in (4.3) gives rise to the term
$$c_1(-1)^{m_i\cN_i}(\cN_i!)^{m_i}
\frac{x^{(1)}(t_i)^{m_i\cN_i}} { x(t_i)^{m_i(\cN_i+1)}},$$ which
is the only term in $P^{m_i}(x^{-1})$ with this high or higher
power of $x(t_i)$ in the denominator. Hence taking
$$x(t)=x_s(t)=e^{2\pi i(t-t_i)}-1+s, \hspace{1mm} s \in (-1,0),$$
in (4.4) and noting that $\varphi(x_s)=\prod_j x_s(t_j)^{n_j}$, we
obtain
$$h_{0 \infty}^{n-m_i}(x_s)=
s^{n_i-m_i(\cN_i+1)}\prod_{j \neq i }x_s(t_j)^{n_j}
\left(c_2+sg(s)\right), \leqno{(4.5)}$$ where $c_2 \neq 0$ is a
constant, $g(s)$ is a rational function in $s$ which is bounded on
the interval $(-1,0)$. In order that the limit of right hand side
of (4.5) exist as $s \to 0^{-}$ we must have $m_i(\cN_i+1) \le
n_i$.

The proof is finished if we can show that for any $\sigma \in
H^0(L_{k}\bP_1, \lp)$ there exists $\sigma^{\prime} \in
H^0(L_{k}\bP_1, \lp)$ such that the corresponding
$\cN_i^{\prime}=\cN_i$, $m_i^{\prime}=m_i$, and
$\hp(\sigma^{\prime})$ contains the monomial (4.3).

We write $P$ in the form
$$P(x)=x^{(\cN_i)}(t_i)^{m_i} f_1(x)+f_2(x),
\hspace{1mm} x \in L_{k}\bC,$$ where $f_1(x), f_2(x)$ are
polynomials in $x^{(\nu)}(t_j)$, $f_1 \neq 0$ is independent of
$x^{(\cN_i)}(t_i)$ and the degree of $f_2$ in $x^{(\cN_i)}(t_i)$
is strictly less than $m_i$. Let $y \in L_{k}\bC^{\ast}$, and
$\gamma \in L_{k}PSL(2, \bC)$ be the map $x \mapsto y(x+1)$, $x
\in L_{k}\bP_1$. Let $P^{\prime}$ denote
$\hp(\gamma^{\ast}\sigma)$. From Propositions 2.7 and 2.5 we
obtain that $\gamma^{\ast}\sigma_{\infty}=c_3 \sigma_{\infty}$
(with a fixed isomorphism $\gamma^{\ast}\lp \cong \lp$), where
$c_3 \neq 0$ is a constant. Therefore
\begin{align*}P^{\prime}(x)=c_3\frac{\gamma^{\ast}\sigma}
{\gamma^{\ast}\sigma_{\infty}}(x)=c_3\gamma^{\ast}
P(x)=c_3P\left(y(x+1)\right), \hspace{1mm} x \in L_{k}\bC.
\tag{4.6}\end{align*} Computing the right hand side of (4.6) by
the product rule, we find that it contains the monomial
$$c_3\, y(t_i)^{m_i}f_1(y)\left(x^{(\cN_i)}(t_i)\right)^{m_i}.$$
Since $f_1 \neq 0$, we can find $y$ such that $f_1(y)\neq 0$; then
we can choose $\sigma^{\prime}=\gamma^{\ast}\sigma$.
\end{proof}

\noindent {\it Proof of Theorem 1.2(1).} The lower bound is given
in Proposition \ref{pb}. For the case of $C^{\infty}$ loops
Propositions \ref{poly} and \ref{up} imply that $\R)$ consists of
polynomials of degree $\le \ord(\lp)$ in $x^{(\nu)}(t_j)$, $0 \le
\nu \le n_j-1$, $1 \le j \le r$, therefore $\dim
H^0(L_{\infty}\bP_1, \lp)< \infty$. Since in general
$L_{\infty}\bP_1$ is continuously embedded into $L\bP_1$ with a
dense image, so that the restriction map $\s \to
H^0(L_{\infty}\bP_1, \lp)$ is monomorphic, we conclude that $\dim
\s < \infty$. \qed
\medskip

An immediate application of Proposition \ref{up} is to identify
all holomorphic sections of a ``generic'' bundle of the type
considered in Theorem 1.2(1).

\begin{cor}
If $n_1=\cdots = n_r=1$ in (4.1), then $\dim \s=2^r$.
\end{cor}
\begin{proof}
Proposition \ref{up} gives that $\R$ only contains polynomials in
evaluation maps $E_{t_1}, \cdots, E_{t_r}$, of degree $\le 1$ in
each variable, and all such polynomials are indeed in $\R$ by
Proposition \ref{pb}.
\end{proof}

There are bundles as in Theorem 1.2(1) which have holomorphic
sections other than those identified in Proposition \ref{pb}. We
shall show this by an explicit construction.

\begin{prop} \label{sec}
Let $h \in \co(L\bC)$ and $\varphi \in \h$. Then $h \in \R$ if and
only if there exists a family $\{f_{a} \in \co(L\bC): a \in \bC\}$
such that
$$h(x+a)=\varphi(x)f_{a}(x^{-1}), \hspace{1mm} x \in
L\bC^{\ast}.\leqno{(4.7)}$$
\end{prop}
\begin{proof}
With $\sigma_a$ of Proposition 2.7 (4.7) is equivalent to
$$h(x)\sigma_{\infty}(x)=f_a\left((x-a)^{-1}\right)\sigma_a(x),
\hspace{1mm} x \in LU_a \cap LU_{\infty}, \hspace{1mm} a \in
\bC,$$ which defines a section $\sigma \in \s$ such that
$\hp(\sigma)=h$ and vice versa.
\end{proof}

From now on, until the last paragraph, we shall work with the
space of $C^k$ loops. Let $h \in \co(L_k\bC^{\ast})$ be a rational
function of finitely many $x^{(\nu)}(t_j)$, $0\le \nu \le k$,
$1\le j\le r$. Letting each $t_j$ vary in $S^1$, $h$ induces a
function $\chi \in C(L_k\bC^{\ast} \times (S^1)^r)$. If $h$ does
not depend on the highest derivative $x^{(k)}(t_i)$ for some $i$
(this is automatically satisfied if $k=\infty$), then $\chi$ is
differentiable in $t_i$, and we define
$$h^{t_i}(x)=\frac{\partial}{\partial t_i}\chi(x, t_1, \cdots,
t_r).$$ Thus $h^{t_i} \in \co(L_k\bC^{\ast})$ is also a rational
function in $x^{(\nu)}(t_j)$, $0\le \nu \le k$, $1\le j \le r$. If
$h$ was a polynomial, so will be $h^{t_i}$. The linear map $T^i: h
\mapsto h^{t_i}$ satisfies $T^i(h_1h_2)=T^i(h_1)h_2+h_1T^i(h_2)$,
$T^i(h(x^{-1}))=(T^ih)(x^{-1})$, and $T^i(h(x+a))=(T^ih)(x+a)$ if
$h \in \co(L_k\bC)$.

\begin{prop} \label{der}
Let $\varphi_i=\varphi E_{t_i}$, and assume that $h \in \R$ does
not depend on $x^{(k)}(t_i)$ for some $i$. Then $h^{t_i} \in
\mathrm{R}(\fH_{\varphi_i})$.
\end{prop}
\begin{proof}
If $f_a$ is defined as in Proposition 4.5, Proposition 3.2 implies
$f_a \in \cL^n(L\bC)$. Since $h$ and $\varphi$ depend only on
$x^{(\nu)}(t_j)$, so does $f_a$: it is in fact a polynomial in
$x^{(\nu)}(t_j)$ ($\nu <k$ if $j=i$, $\nu \le k$ for all other
$j$). Applying $T^i$ to (4.7) we obtain
\begin{align*}
h^{t_i}(x+a)=T^i\big(h(x+a)\big)
&=n_i\varphi(x)x(t_i)^{-1}x^{(1)}(t_i)f_{a}(x^{-1})
+\varphi(x)f_{a}^{t_i}(x^{-1})\\&=\varphi_i(x)
\tilde{f}_a(x^{-1}),
\end{align*} where
$\tilde{f}_a(x)=-n_ix^{(1)}(t_i)f_{a}(x)+x(t_i)f_{a}^{t_i}(x)$, so
that $h^{t_i} \in \mathrm{R}(\fH_{\varphi_i})$ by Proposition
\ref{sec}.
\end{proof}

Finally we shall discuss the order $N_i$ of the highest derivative
$x^{(\nu)}(t_i)$ that $\hs$ can depend on, in the case of the
$C^k$ loop space $L_k\bP_1$. Recall from Propositions 4.2, 4.3
that $N_i \le \min(k,n_i-1)$. It turns out that this estimate is
sharp:

\begin{thm} \label{exc}
There exists a section of $\lp$ for which $N_i=\min(k, n_i-1)$.
\end{thm}
\begin{proof}
With $\mu_i=\min(k,n_i-1)$ and $0 \le \nu \le \mu_i$ let
$\varphi_{\nu}=\varphi E_{t_i}^{\nu-\mu_i}$. Proposition 4.1 gives
that $E_{t_i} \in \mathrm{R}(\fH_{\varphi_0})$. Repeatedly
applying Proposition 4.6 we obtain that the functions
$h_{\nu}(x)=x^{(\nu)}(t_i)$ are in
$\mathrm{R}(\fH_{\varphi_{\nu}})$ for $\nu \le \mu_i$. Thus
$h_{\mu_i}(x)=x^{(\mu_i)}(t_i)$ is in the range of $\hp$, and
$N_i=\mu_i$ if $\sigma=\hp^{-1}h_{\mu_i}$.
\end{proof}

A similar reasoning would apply to the space of $W^{k,p}$ loops,
where the largest value for $N_i$ turns out to be $\min(k-1,
n_i-1)$. A remarkable consequence of this is that while for
generic bundles $\lp$, as we have seen in Corollary 4.4, $\dim
H^0(L\bP_1,\lp)$ does not vary with the regularity of loops, when
at least one $n_i >1$, this dimension will depend on the
regularity class $C^k$, $W^{k,p}$ considered.

\end{document}